\documentclass[a4paper,12pt]{amsart}
\usepackage{hyperref,fancyhdr,mathrsfs,amsmath,amscd,amsthm,amsfonts,latexsym,amssymb}

\DeclareMathOperator{\dif}{d}

\renewcommand{\H}{\mathscr{H}}

\newcommand{\F}{\mathscr{F}}
\newcommand{\Fa}{\mathcal{F}}
\newcommand{\J}{\mathcal{J}}
\newcommand{\tG}{\mathscr{G}}

\def \G{\Gamma}

\def \phi{\varphi}
\def \Phi{\varPhi}
\def \p{\pi}

\def \t{\tau}

\def \R{\mathbb{R}}

\def \C{\mathbb{C}\,}

\def\widecheckg{g^{\hspace*{-2.5pt}\vbox to 5pt{\hbox to
0pt{\LARGE$\check{}$}}}\hspace*{2pt}}

\def\widecheckl{\lambda^{\hspace*{-3.5pt}\vbox to 8pt{\hbox to
0pt{\LARGE$\check{}$}}}\hspace*{2pt}}

\begin{document}

\title{On a class of twistorial maps}
\author{Radu~Pantilie}
\thanks{The author gratefully acknowledges that this work was partially supported by a
CEx grant no.\ 2-CEx 06-11-22/25.07.2006.}
\email{\href{mailto:Radu.Pantilie@imar.ro}{Radu.Pantilie@imar.ro}}
\address{R.~Pantilie, Institute of Mathematics ``Simion Stoilow'' of the Romanian Academy,
P.O. BOX 1-764, RO-014700 Bucharest, Romania}

\newtheorem{thm}{Theorem}[section]
\newtheorem{lem}[thm]{Lemma}
\newtheorem{cor}[thm]{Corollary}
\newtheorem{prop}[thm]{Proposition}

\theoremstyle{definition}

\newtheorem{defn}[thm]{Definition}
\newtheorem{rem}[thm]{Remark}
\newtheorem{exm}[thm]{Example}

\numberwithin{equation}{section}

\maketitle
\thispagestyle{empty}
\vspace{-0.5cm}
\section*{Abstract}
\begin{quote}
{\footnotesize We show that a natural class of twistorial maps gives a pattern for apparently different
geometric maps, such as, $(1,1)$-geodesic immersions
from $(1,2)$-symplectic almost Hermitian manifolds and pseudo horizontally conformal
submersions with totally geodesic fibres for which the associated almost CR-structure is integrable.
Along the way, we construct for each constant curvature Riemannian manifold $(M,g)$\,, of dimension
$m$\,, a family of twistor spaces $\bigl\{Z_r(M)\bigr\}_{1\leq r<\tfrac12m}$ such that $Z_r(M)$
parametrizes naturally the set of pairs $(P,J)$\,, where $P$ is a totally geodesic submanifold of $(M,g)$\,,
of codimension $2r$\,, and $J$ is an orthogonal complex structure on the normal bundle of $P$ which
is parallel with respect to the normal connection.}
\end{quote}

\section*{Introduction}

\indent
In the complex-analytic category, the \emph{twistor space} of a manifold $M$, endowed with a twistorial structure,
parametrizes the set of certain submanifolds --\:\emph{the twistors}\:-- of $M$. For example (see \cite{PanWoo-sd}
and the references therein),
the twistor space of a three-dimensional complex Einstein--Weyl space $(M^3,c,D)$ is formed of the (maximal)
degenerate surfaces of $(M^3,c)$ which are totally geodesic with respect to $D$\,. Also, the twistor space of
a four-dimensional anti-self-dual complex-conformal manifold $(M^4,c)$ is formed of the self-dual surfaces
of $(M^4,c)$ (similar comments apply, for example, to the complex-quaternionic
manifolds of dimension at least eight). Further, the space of (unparametrized) isotropic geodesics of a
complex-conformal manifold is, in a natural way, a twistor space \cite{LeB-nullgeod}\,.\\
\indent
In the smooth category, the definition of almost twistorial structure is slightly different \cite{LouPan-II}\,;
it follows that, \emph{in the smooth category, the twistors are certain submanifolds for which the normal bundle
is endowed with a (linear) CR-structure}.
These submanifolds may well be just points. For example, the twistor space $Z$ of a three-dimensional
conformal manifold $(M^3,c)$ is a five-dimensional CR-manifold formed of the orthogonal nontrivial CR-structures
on $(M^3,c)$ \cite{LeB-CR} (assuming $(M^3,c)$ real-analytic, $Z$ is a real hypersurface, endowed with the induced
CR-structure, in the space of isotropic geodesics of the germ-unique complexification of $(M^3,c)$\,).
Also, it is well-known (see \cite{LouPan-II} and the references therein) that the twistor space
of a four-dimensional anti-self-dual conformal manifold $(M^4,c)$ is a complex manifold, of complex dimension three,
formed of the positive orthogonal complex
structures on $(M^4,c)$ (similar comments apply, for example, to the quaternionic manifolds of dimension at least
eight).\\
\indent
On the other hand (see \cite{LouPan-II} and the references therein),
the twistor space $Z$ of a three-dimensional Einstein--Weyl space $(M^3,c,D)$ is, locally, a
complex surface formed of the pairs $(\gamma,J)$\,, where $\gamma$ is a geodesic of $D$ and $J$ is an
orthogonal complex structure
on the normal bundle of $\gamma$ (obviously, if $M^3$ is oriented then $Z$ is just the space of oriented geodesics
of $D$). In \S\ref{section:f}\,, below, we generalize this example by constructing for each constant curvature
Riemannian manifold $(M,g)$\,, of dimension $m\geq3$\,, a family of twistor spaces
$\bigl\{Z_r(M)\bigr\}_{1\leq r<\tfrac12m}$ such that $Z_r(M)$ is, locally, a complex manifold, of complex dimension
$r(2m-3r+1)/2$\,, formed of the pairs $(P,J)$\,, where $P$ is a totally geodesic submanifold of $(M,g)$\,,
of codimension $2r$\,, and $J$ is an orthogonal complex structure on the normal bundle of $P$ which
is parallel with respect to the normal connection (the particular case $r=1$ is due to \cite{BaiWoo2}\,).
Moreover, we prove that, in dimension at least four, the constant curvature Riemannian manifolds give
all the Weyl spaces for which this construction works (Theorem \ref{thm:tcc}\,).\\
\indent
A map $\phi:M\to N$ between manifolds endowed with twistorial structures is \emph{twistorial} if it maps
consistently (some of the) twistors on $M$ to twistors on $N$ (see \S\ref{section:twiststrmap_class} for a definition suitable
for this paper and \cite{LouPan-II} for a more general definition; cf.\ \cite{PanWoo-sd}\,).\\
\indent
In this paper, we show that apparently different geometric maps are examples of such twistorial maps:\\ 
\indent
\:\:$\bullet$ in \S\ref{section:imm}\,, we study twistorial immersions between even dimensional oriented
Weyl spaces endowed with the
\emph{associated nonintegrable almost twistorial structures} (see Example \ref{exm:flatnonint}\,;
cf.\ \cite{EelSal}\,),\\
\indent
\:\:$\bullet$ in \S\ref{section:(1,1)submfds}\,, we study $(1,1)$-geodesic immersions
from $(1,2)$-symplectic almost Hermitian manifolds,\\
\indent
\:\:$\bullet$ in \S\ref{section:f}\,, we study pseudo horizontally conformal submersions with totally geodesic
fibres for which the associated almost CR-structure is integrable.\\
\indent
In \S\ref{section:int}\,, we prove a powerful integrability result (Theorem \ref{thm:O'BrRaw}\,;
cf.\ \cite{O'BrRaw}\,) which can be applied to all of the examples of almost twistorial structures known to us.
Here, we use this result to give the necessary and sufficient conditions
for the integrability of several almost CR-structures and almost $f$-structures (see, for example,
Theorems \ref{thm:LeBr-O'BrRaw} and \ref{thm:tcc}\,), related to the twistorial maps we consider.\\
\indent
See \cite{PanWoo-sd}\,, \cite{LouPan} and \cite{LouPan-II} for more information on almost
twistorial structures and twistorial maps.\\
\indent
I am grateful to John~C.~Wood, and Eric~Loubeau for useful comments, and to Liviu~Ornea for useful discussions.

\section{An integrability result} \label{section:int}

\indent
Let $F_j$ be a complex submanifold of the Grassmannian manifold ${\rm Gr}_{r_j}(m_j,\C)$\,, $1\leq r_j\leq m_j$\,, $(j=1,2)$\,.
Suppose that there exists a complex Lie subgroup $G_j$ of ${\rm GL}(m_j,\C)$ whose canonical action on ${\rm Gr}_{r_j}(m_j,\C)$
induce a tranzitive action on $F_j$\,; thus, $F_j=G_j/H_j$\,, as complex manifolds, where $H_j$ is the isotropy group of
$G_j$ at some point of $F_j$\,, $(j=1,2)$\,.\\
\indent
Let $(P_j,M,G_j)$ be a (smooth) principal bundle endowed with a connection $\nabla_j$\,, $(j=1,2)$\,, where $M$ is a (smooth connected) manifold,
$\dim M=m_1$\,. We suppose that $(P_1,M,G_1)$ is a subbundle of the complex frame bundle of $T^{\C}\!M$.
Denote $Q_j=P_j\times_{G_j}F_j$ and
let $\H^j\subseteq TQ_j$ be the connection induced by $\nabla_j$ on $Q_j$\,; note that, $Q_j=P/H_j$\,, $(j=1,2)$\,.\\
\indent
Denote $Q=\iota^*(Q_1\times Q_2)$ where $\iota:M\to M\times M$ is defined by $\iota(x)=(x,x)$\,, for any $x\in M$,
and let $\p:Q\to M$ be the projection. Obviously, ${\rm ker}\dif\!\p$ is a complex vector bundle.\\
\indent
The connections $\nabla_1$ and $\nabla_2$ induce a connection $\H\,(\subseteq TQ)$ on $Q$ and
let $\tG_0\subseteq\H^{\C}$ be the subbundle characterised by $\dif\!\p(\tG_0)(p,q)=p$\,, for all
$(p,q)\in Q$\,. Define
\begin{equation*}
\begin{split}
\tG=&\tG_0\oplus({\rm ker}\dif\!\p)^{0,1}\;,\\
\tG'=&\tG_0\oplus({\rm ker}\dif\!\p)^{1,0}\;.
\end{split}
\end{equation*}
\indent
Let $T$ be the torsion of $\nabla_1$ and let $R_j$ be the curvature form of $\nabla_j$\,, $(j=1,2)$\,.

\begin{thm}[cf.\ \cite{O'BrRaw}\,, \cite{EelSal}\,] \label{thm:O'BrRaw}
If\/ $\dim_{\C}\!F_1\geq1$ then $\tG'$ is nonintegrable. Furthermore, the following assertions are equivalent:\\
\indent
{\rm (i)} $\tG$ is integrable.\\
\indent
{\rm (ii)} $T(\Lambda^2 p)\subseteq p$\,, $R_1(\Lambda^2p)(p)\subseteq p$\,, $R_2(\Lambda^2p)(q)\subseteq q$\,,
for all $(p,q)\in Q$\,.
\end{thm}
\begin{proof}
Let $G=G_1\times G_2$ and let $\mathfrak{g}$ and $\mathfrak{g}_j$ be the Lie algebras of $G$ and $G_j$\,, respectively,
$(j=1,2)$\,. Let $P=\iota^*(P_1\times P_2)$ and denote, also, by $\p$ and $\H$ the projection $\p:P\to M$ and
the connection $\H\subseteq TP$ induced by $\nabla_1$ and $\nabla_2$ on $P$. Note that,
$({\rm ker}\dif\!\p)^{\C}=P\times(\mathfrak{g}\oplus\overline{\mathfrak{g}})$\,.\\
\indent
Let $H_1$ and $H_2$ be the isotropy groups of $G_1$ and $G_2$ at some points $p_0$ and $q_0$ of $F_1$ and $F_2$\,, respectively.
Let $H=H_1\times H_2$ and denote by $\mathfrak{h}$ its Lie algebra.\\
\indent
Let $\tG_{0,P}\subseteq\H^{\C}$ be characterised by $(\dif\!\p)_{(u,v)}(\tG_{0,P})=u(p_0)$\,, for all $(u,v)\in P$.
Clearly, $\dif\!\psi(\tG_{0,P})=\tG_0$ and $(\dif\!\psi)^{-1}(\tG_0)=\tG_{0,P}$ where $\psi:P\to Q$ is the projection.\\
\indent
Define
\begin{equation}
\begin{split}
\tG_P&=\tG_{0,P}\oplus P\times\mathfrak{h}\oplus P\times\overline{\mathfrak{g}}\;,\\
\tG_P'&=\tG_{0,P}\oplus P\times\overline{\mathfrak{h}}\oplus P\times\mathfrak{g}\;.
\end{split}
\end{equation}
Obviously, $\dif\!\psi(\tG_P)=\tG$\,, $\dif\!\psi(\tG_P')=\tG'$\,, $(\dif\!\psi)^{-1}(\tG)=\tG_P$\,, $(\dif\!\psi)^{-1}(\tG)=\tG_P'$\,.
Therefore, $\tG$ is integrable if and only if $\tG_P$ is integrable. Similarly, $\tG'$ is nonintegrable if and only if
$\tG_P'$ is nonintegrable.\\
\indent
For each $\xi\in\C^{\!m_1}$ let $B(\xi)$ be the horizontal (complex) vector field on $P$ characterised by $(\dif\!\p)_{(u,v)}(B(\xi))=u(\xi)$
for any $(u,v)\in P$. Obviously, the map $P\times p_0\to\tG_{0,P}$\,, $((u,v),\xi)\mapsto B(\xi)_{(u,v)}$\,, for $(u,v)\in P$ and $\xi\in p_0$\,,
is an isomorphism of vector bundles. Also, if $A\in\mathfrak{g}_1$ and $\xi\in\C^{\!m_1}$
then $[A,B(\xi)]=B(A\,\xi)$ and $[\overline{A},B(\xi)]=0$ (cf.\ \cite[Chapter III]{KoNo}\,)\,.\\
\indent
If $\dim_{\C}\!F_1\geq1$ then for any $A\in\mathfrak{g}_1\!\setminus\!\mathfrak{h}_1$ and $\xi\in p_0$ we have $A\,\xi\notin p_0$\,.
Hence, $[A,B(\xi)]=B(A\,\xi)$ is nowhere tangent to $\tG_{0,P}$\,. Therefore $\tG_P'$ and $\tG'$ are nonintegrable.\\
\indent
The equivalence (i)$\iff$(ii) follows straightforwardly from Cartan's structural equations (cf.\ \cite{Pan-Deva}\,).
\end{proof}

\indent
Let $p$ be a section of $Q_1$\,; we shall denote by the same letter $p$ the corresponding complex vector subbundle of $T^{\C}\!M$.
Then the map $Q_2\hookrightarrow Q$\,, $q\mapsto(p_{\p_2(q)},q)$\,, for any $q\in Q_2$\,, is an embedding, where
$\p_2:Q_2\to M$ is the projection. Denote by the same symbol $Q_2$ the image of this embedding. Then
$\tG^p=\tG\cap TQ_2$ is a subbundle of $T^{\C}\!Q_2$\,.\\
\indent
Note that, $\tG^p$ does not depend of $\nabla_1$. In fact, we could define $\tG^p$ as follows. Firstly, let
$\tG^p_0$ be the subbundle of $T^{\C}\!Q_2$ which is horizontal, with respect to the connection induced by $\nabla_2$ on $Q_2$\,,
and such that $\dif\!\p_2(\tG^p_0)=p$. Then we have $\tG^p=\tG^p_0\oplus({\rm ker}\dif\!\p_2)^{0,1}$\,.\\
\indent
{}From Theorem \ref{thm:O'BrRaw} we easily obtain the following result.

\begin{cor}[cf.\ \cite{KosMal}\,] \label{cor:O'BrRaw}
The following assertions are equivalent:\\
\indent
{\rm (i)} $\tG^p$ is integrable.\\
\indent
{\rm (ii)} $p$ is integrable and $R_2(\Lambda^2_xp)(q)\subseteq q$ for any $x\in M$ and $q\in\p_2{}^{-1}(x)$\,.
\end{cor}

\begin{rem}
Theorem \ref{thm:O'BrRaw} and Corollary \ref{cor:O'BrRaw} can be easily generalized to the case when $Q_2$ is a fibre bundle
for which the typical fibre is a complex manifold and the structural group is a complex Lie group whose action on the typical fibre
is transitive and holomorphic.
\end{rem}

\section{Almost twistorial structures and twistorial maps} \label{section:twiststrmap_class}

\indent
An \emph{almost CR-structure} on a (smooth connected) manifold $M$ is a section $J$ of ${\rm End}(\H)$ such that
$J^2=-\,{\rm Id}_{\H}$\,, for some distribution $\H$ on $M$\,; if $\H=TM$ then $J$ is an \emph{almost complex structure}
on $M$. Let $\F$ be the eigenbundle of (the complexification of) $J$ corresponding to $-{\rm i}$\,; we say that $\F$ is
\emph{the complex distribution associated to $J$}\,.
Then $J$ is \emph{integrable} if $\F$ is integrable (that is, for any $X,Y\in\G(\F)$ we have $[X,Y]\in\G(\F)$\,).
A \emph{CR-structure} is an integrable almost CR-structure; a \emph{complex structure} is an integrable
almost complex structure (see \cite[\S2]{LouPan-II}\,).\\
\indent
An \emph{almost $f$-structure} on $M$ is a section $F$ of ${\rm End}(TM)$ such that $F^3+F=0$\,.
Let $\F=T^0M\oplus T^{0,1}M$ where $T^0M$ and $T^{0,1}M$ are the eigenbundles of $F$ corresponding
to $0$ and $-{\rm i}$\,, respectively; we say that $\F$ is \emph{the complex distribution associated to $F$}\,.
Then $F$ is \emph{integrable} if $\F$ is integrable. An \emph{$f$-structure} is an integrable almost $f$-structure
(see \cite[\S2]{LouPan-II}\,).\\
\indent
If $\F$ is the complex distribution associated to a CR-structure or an $f$-structure on a manifold $M$
then, obviously, $\F\cap\overline{\F}$ is (the tangent bundle of) a foliation on $M$.\\
\indent
Let $F$ be an almost $f$-structure on $M$ and let $T^{1,0}M$ and $T^{0,1}M$ be its eigenbundles
corresponding to ${\rm i}$ and $-{\rm i}$\,, respectively. Then $J=F|_{T^{1,0}M\oplus T^{0,1}M}$ is an almost
CR-structure on $M$; we shall call $J$ \emph{the almost CR-structure induced by $F$}\,.
Note that, $F$ is not determined by $J$; also, if $F$ is integrable then $J$ is not
necessarily integrable.\\
\indent
An \emph{(almost) CR-structure on a conformal manifold} $(M,c)$ is an (almost) CR-structure $J$ on $M$
such that $J^*+J=0$\,; obviously, this holds if and only if the complex distribution associated to $J$
is isotropic.\\
\indent
An \emph{(almost) $f$-structure on a conformal manifold} $(M,c)$ is an (almost) $f$-structure $F$ on $M$
such that $F^*+F=0$\,; obviously, this holds if and only if $T^{0,1}M$ is isotropic and
$T^0M=\bigl(T^{1,0}M\oplus T^{0,1}M\bigr)^{\perp}$\,.
Therefore an almost $f$-structure on a conformal manifold is determined
by its eigenbundle corresponding to ${\rm i}$ (or ${-\rm i}\,)$. Equivalently, if we denote by $J$ the
almost CR-structure whose eigenbundle corresponding to ${\rm i}$ is $T^{1,0}M$ then $F\longleftrightarrow J$
establishes a bijective correspondence (which depends on $c$\,) between almost $f$-structures on $(M,c)$
and almost CR-structures on $(M,c)$\,.

\begin{defn}
A map $\phi:(M,F^M)\to(N,F^N)$\,, between manifolds endowed with almost $f$-structures (or, almost
CR-structures), is \emph{holomorphic} if $\dif\!\phi(\F^M)\subseteq\F^N$\,, where $\F^M$ and $\F^N$
are the complex distributions associated to $F^M$ and $F^N$, respectively.
\end{defn}

\begin{rem}
An almost $f$-structure $F$ on $M$ is integrable if and only if for any $x\in M$ there exists an open
neighbourhood $U\ni x$ and a holomorphic submersion $\phi$ from $(U,F|_U)$ onto some complex manifold
$(N,J)$ such that ${\rm ker}\dif\!\phi=T^0M$ \cite{Pan-Deva}\,; we say that the $f$-structure $F|_U$ is
\emph{defined by $\phi$}\,. A \emph{simple $f$-structure} is an $f$-structure
(globally) defined by a holomorphic submersion with connected fibres.
\end{rem}

\indent
We end this section with the definitions of almost twistorial structure and twistorial map
suitable for the purpose of this paper; more general definitions are given in \cite{LouPan-II}
(cf.\ \cite{PanWoo-sd}\,).

\begin{defn}
An \emph{almost twistorial structure} on a manifold $M$ is a quadruple $\t=(Q,M,\p,\J)$\,, where $\p:Q\to M$ is a
locally trivial fibre space and $\J$ is an almost CR-structure or an almost $f$-structure on $Q$
which induces almost complex structures on each fibre of $\p$\,.
We say that $\t$ is \emph{integrable} if $\J$ is integrable; a \emph{twistorial structure}
is an integrable almost twistorial structure. Suppose that $\t$ is a twistorial structure such that
there exists a surjective submersion $\phi:Q\to Z$ whose fibres are the leaves of $\F\cap\overline{\F}$\,,
where $\F$ is the complex distribution associated to $\J$. Then $Z$, endowed with the CR-structure
$\dif\!\phi(\F)$\,, is the \emph{twistor space} of $\t$\,.
\end{defn}

\begin{defn}
Let $\phi:M\to N$ be a map between manifolds endowed with the almost twistorial structures $\t_M=(Q_M,M,\p_M,\J^M)$ and
$\t_N=(Q_N,N,\p_N,\J^N)$\,. Suppose that there exists a section $p$ of $Q_M$ and a map $\Phi:p(M)\to Q_N$ such that
$\p_N\circ\Phi=\phi\circ\p_M|_{p(M)}$ and the tangent bundle of $p(M)$ is preserved by $\J^M$; denote by
$\J^p$ the restriction of $\J^M$ to the tangent bundle of $p(M)$\,.
We shall say that $\phi:(M,\t_M)\to(N,\t_N)$ is a \emph{twistorial map (with respect to $\Phi$)},
if $\Phi:(p(M),\J^p)\to(Q_N,\J^N)$ is holomorphic; that is, $\dif\!\Phi(\F^p)\subseteq\F^N$
where $\F^p$ and $\F^N$ are the complex distributions associated to $\J^p$ and $\J^N$, respectively.
\end{defn}

\section{Twistorial immersions between Weyl spaces} \label{section:imm}

\indent
We start this section with two related examples of almost twistorial structures.

\begin{exm} \label{exm:flatnonint}
Let $(M,c,D)$ be an oriented, even-dimensional Weyl space and let $\p:Q\to M$ be the bundle of positive maximal isotropic spaces
on $(M,c)$ (the positive maximal isotropic spaces on $(M,c)$ are the eigenspaces, corresponding to $-{\rm i}$\,, of the positive orthogonal
complex structures on $(M,c)$\,).
As ${\rm ker}\dif\!\p$ is a complex vector bundle, we have an isomorphism of complex vector bundles
$({\rm ker}\dif\!\p)^{\C}=({\rm ker}\dif\!\p)^{1,0}\oplus({\rm ker}\dif\!\p)^{0,1}$.
Let $\H\subseteq TQ$ be the connection induced by $D$ on $Q$\,. Let $\tG^0\subseteq\H^{\C}$ be the complex vector subbundle characterised by
$\dif\!\p(\tG^0_q)=q$\,, for any $q\in Q$\,, and define
\begin{equation*}
\begin{split}
\tG=&\tG_0\oplus({\rm ker}\dif\!\p)^{0,1}\;,\\
\tG'=&\tG_0\oplus({\rm ker}\dif\!\p)^{1,0}\;.
\end{split}
\end{equation*}
\indent
Let $\J$ and $\J'$ be the almost complex structures whose eigenbundles corresponding to $-{\rm i}$ are $\tG$ and $\tG'$,
respectively.\\
\indent
Obviously, if $\dim M=2$ then $Q=M$ and $\J=\J'$ is the positive Hermitian structure of $(M^2,c)$\,.\\
\indent
Note that, $\J$ does not depend of $D$ whilst if $\dim M\geq4$ then $\J'$ determines $D$ (that is, if $D_1$ is another
Weyl connection on $(M,c)$ which induces $\J'$ then $D=D_1$\,; this follows from \cite[Proposition 2.6]{LouPan-II}\,).\\
\indent
If $\dim M\geq4$ then $\J'$ is nonintegrable (that is, always not integrable) whilst if $\dim M=4$ then $\J$ is
integrable if and only if $(M^4,c)$ is anti-self-dual and if $\dim M\geq6$ then $\J$ is integrable if and only
if $(M,c)$ is flat; these well-known results (see \cite[\S4]{EelSal}\,, \cite[\S5]{O'BrRaw}\,, \cite[\S3]{Pan-Deva}\,)
follow from Theorem \ref{thm:O'BrRaw}\,.\\
\indent
Obviously, $(Q,M,\p,\J)$ and $(Q,M,\p,\J')$ are almost twistorial structures on $M$; we shall call
$(Q,M,\p,\J')$ the \emph{nonintegrable almost twistorial structure associated to $(M,c,D)$}\,.
\end{exm}

\indent
Let $(M,c_M,D^M)$ and $(N,c_N,D^N)$ be even-dimensional oriented Weyl spaces and let $\t'_M=(Q_M,M,\p_M,\J'_M)$ and $\t'_N=(Q_N,N,\p_N,\J'_N)$
be the associated nonintegrable almost twistorial structures.\\
\indent
Suppose that $\phi:M\hookrightarrow N$ is an injective immersion. Then orient $(TM)^{\perp}$ such that the isomorphism
$TN|_M=TM\oplus(TM)^{\perp}$ be orientation preserving and let $\p:Q\to M$ be the bundle of positive maximal isotropic spaces
on $\bigl((TM)^{\perp},c_N|_{(TM)^{\perp}}\bigr)$\,.\\
\indent
If $p$ is a (local) section of $Q_M$ then we shall denote by $J^p$ the almost Hermitian structure on $(M,c_M)$ such that
$p$ is the eigenbundle of $J^p$ corresponding to $-{\rm i}$\,; similarly, for $Q_N$\,.
Standard arguments show that the following assertions are equivalent:\\
\indent
\quad(i) $p:(M,J^p)\to(Q_M,\J'_M)$ is holomorphic.\\
\indent
\quad(ii) $D^M_{\overline{X}}Y$ is a section of $p$ for any sections $X$ and $Y$ of $p$\,.\\
\indent
\quad(iii) $D^M_{J^pX}J^p=-J^pD^M_{\,X}\!J^p$, for any $X\in TM$.\\
\indent
\quad(iv) $(\dif^{D^M}\!\omega_M)^{(1,2)}=0$\,, where $\omega_M$ is the K\"ahler form of $(M,c_M,J^p)$ (defined by
$\omega_M(X,Y)=c_M(J^pX,Y)$\,, for any $X,Y\in TM$).\\
Furthermore, if assertion (i)\,, (ii) or (iii) holds then $D^M$ is the Weyl connection of $(M,c_M,J^p)$ (see \cite[Remark 3.3]{LouPan}\,).\\
\indent
We shall denote by $\J^p$ the almost complex structure on $Q$ whose eigenbundle corresponding to ${\rm i}$ is constructed,
similarly to $\tG^p$ of Corollary \ref{cor:O'BrRaw}\,, by using the connection induced by $\Pi\circ D^N$ on $Q$ and the
complex vector subbundle $\overline{p}$ of $T^{\C}\!M$, where $\Pi:TN|_M\to(TM)^{\perp}$ is the orthogonal projection.\\
\indent
Let $L$ be the line bundle of $(N,c_N)$\,. We define a section $A$ of the bundle
$(L|_M)^2\otimes\Lambda^2T^*M\otimes\Lambda^2\bigl((TM)^{\perp}\bigr)^*$ by
$$A(X,Y,U,V)=\sum_ac_N(D^N_{\,X}Z_a,U)c_N(D^N_{\,Y}Z_a,V)-c_N(D^N_{\,Y}Z_a,U)c_N(D^N_{\,X}Z_a,V)\;,$$
for any $x\in M$ and $X,Y\in T_xM$, $U,V\in (T_xM)^{\perp}$\,, where $\{Z_a\}$ is any conformal local frame on $(M,c_M)$
defined on some open neighbourhood of $x$\,. It is easy to see that $A$ does not depend of $D^N$. Furthermore,
if $M$ is an umbilical submanifold of $(N,c_N)$
then $A=0$\,.

\begin{cor}[cf.\ \cite{SimSve}\,] \label{cor:1}
The almost complex structure $\J^p$ does not depend of the Weyl connection $D^N$. Moreover, the following assertions
are equivalent:\\
\indent
{\rm (i)} $\J^p$ is integrable.\\
\indent
{\rm (ii)} $J^p$ is integrable and $(W+A)(\Lambda^2_xp,\Lambda^2q)=0$ for any $x\in M$ and $q\in Q_x$\,, where $W$ is the Weyl tensor
of $(N,c_N)$\,.
\end{cor}
\begin{proof}
A straightforward calculation gives the following relation, essentially due to Ricci (see \cite[1.72(e)]{Bes}\,),
\begin{equation} \label{e:Ricci}
c_N(R^N\!(X,Y)U,V)=c_N(R^{\Pi}\!(X,Y)U,V)+A(X,Y,U,V)\;,
\end{equation}
for any $X,Y\in TM$ and $U,V\in (TM)^{\perp}$\,, where $R^N$ and $R^{\Pi}$ are the curvature forms of $D^N$ and $\Pi\circ D^N$, respectively.\\
\indent
Also, we have (see \cite{Cal-F}\,)
\begin{equation} \label{e:Weyl}
c_N(R^N\!(X,Y)U,V)=-W(X,Y,U,V)+F^N(X,Y)c_N(U,V)\;,
\end{equation}
for any $X,Y\in TM$ and $U,V\in (TM)^{\perp}$\,, where $W$ is the Weyl tensor of $(N,c_N)$ and $F^N$ is the curvature form of the
connection induced by $D^N$ on $L$\,.\\
\indent
The proof now follows quickly from Corollary \ref{cor:O'BrRaw}\,.
\end{proof}

\begin{rem}
In Corollary \ref{cor:1}\,, if $\dim M=2$ then assertion (ii) is automatically satisfied whilst if ${\rm codim}M=2$ then
the second part of assertion (ii) is automatically satisfied.
\end{rem}

\indent
Let $p$ be a section of $Q_M$ which is isotropic with respect to $c_N$. Then for any map $\Phi:p(M)\to Q_N$ such that
$\p_N\circ\Phi=\phi\circ\p_M|_{p(M)}$ there exists a unique section $q$ of $Q$ such that $\Phi\circ p=p\oplus q$\,.\\
\indent
The following result reduces to \cite[Theorem 5.3]{EelSal}\,, when $\dim M=2$\,, $\dim N=4$
(see \cite[Proposition 5.2]{LouPan-II}\,).

\begin{prop}[cf.\ \cite{SimSve}\,] \label{prop:2}
Let $\Phi$ be given by the sections $p$ and $q$ of $Q_M$ and $Q$\,, respectively, with $p$ isotropic with respect to $c_N$\,.
Then the following assertions are equivalent:\\
\indent
{\rm (i)} $\phi:(M,\t'_M)\to(N,\t'_N)$ is twistorial, with respect to $\Phi$\,.\\
\indent
{\rm (ii)} $\phi$ is $(1,1)$-geodesic with respect to $J^p$ and, $p:(M,J^p)\to(P_M,\J'_M)$ and $q:(M,J^p)\to(Q,\J^p)$ are holomorphic.
\end{prop}
\begin{proof}
Assertion (i) holds if and only if $p(M)$ is an almost complex submanifold of $(Q_M,\J'_M)$ and
$\Phi:\bigl(p(M),\J'_M|_{p(M)}\bigr)\to(Q_N,\J'_N)$ is holomorphic. It is clear that $p(M)$ is an almost complex submanifold of $(Q_M,\J'_M)$
if and only if $p:(M,J^p)\to(P_M,\J'_M)$ is holomorphic. Then $\Phi:\bigl(p(M),\J'_M|_{p(M)}\bigr)\to(Q_N,\J'_N)$ is holomorphic if and only if
$\Phi\circ p:(M,J^p)\to(Q_N,\J'_N)$ is holomorphic. {}From \cite[Proposition 2.6]{LouPan-II} it follows quickly that, assertion (i) is equivalent to
(a) $D^M_{\overline{X}}Y\in\G(p)$\,,
for any $X,Y\in\G(p)$\,, (b) $D^N_{\overline{X}}Y\in\G(p\oplus q)$\,, for any $X,Y\in\G(p)$\,, and
(c) $D^N_{\overline{X}}U\in\G(p\oplus q)$\,, for any $X\in\G(p)\,,\,U\in\G(q)$\,.\\
\indent
Note that, if (b) holds, condition (c) is equivalent to $\Pi(D^N_{\overline{X}}U)\in\G(q)$\,, for any $X\in\G(p)\,,\,U\in\G(q)$\,.
Thus, if (b) holds, condition (c) is equivalent to $q:(M,J^p)\to(Q,\J^p)$ be holomorphic.\\
\indent
Also, if (a) holds, condition (b) is equivalent to $(D\!\dif\!\phi)(\overline{X},Y)\in\G(p\oplus q)$\,, for any $X,Y\in\G(p)$\,.
As $D^M$ and $D^N$ are torsion free, $D\!\dif\!\phi$ is symmetric. It follows quickly that, if (a) holds,
then (b) is equivalent to $(D\!\dif\!\phi)^{(1,1)}=0$\,.\\
\indent
The proposition is proved.
\end{proof}

\begin{rem}
1) If assertion {\rm (i)} or {\rm (ii)} of Proposition \ref{prop:2} holds then $D^M$ is the Weyl connection of
$(M,c_M,J^p)$\,; if, further, $\phi$ is conformal then $D^M$ is equal to the connection induced by $D^N$ on $M$.\\
\indent
2) A  result similar to (but more complicated than) Proposition \ref{prop:2} can be given for twistorial submersions
between Weyl spaces endowed with the nonintegrable almost twistorial structures. It follows again that
such maps are $(1,1)$-geodesic (in particular, harmonic) and, if the codomain is of dimension two, harmonic morphisms.
\end{rem}

\indent
Let $\phi:(M,c_M)\hookrightarrow(N,c_N)$ be a conformal injective immersion. Denote by $Q_M+Q$ the pull-back
by $\iota$ of $Q_M\times Q$\,, where $\iota:M\to M\times M$ is defined by $\iota(x)=(x,x)$\,, for any $x\in M$.
Let $\J$ and $\J'$ be the almost complex structures on $Q_M+Q$ whose eigenbundles corresponding to $-{\rm i}$ are constructed,
similarly to $\tG$ and $\tG'$, respectively, of Theorem \ref{thm:O'BrRaw}\,, by using the connections induced
by $D^M$ and $\Pi\circ D^N$ on $Q_M$ and $Q$\,, respectively.

\begin{prop} \label{prop:umbgeod}
Let $\Phi:Q_M+Q\to Q_N$ be defined by $\Phi(p,q)=p\oplus q$\,, for any $(p,q)\in Q_M+Q$\,.\\
\indent
{\rm (i)} The following assertions are equivalent:\\
\indent
\quad{\rm (i1)} $\Phi:(Q_M+Q,\J)\to(Q_N,\J_N)$ is holomorphic.\\
\indent
\quad{\rm (i2)} $M$ is an umbilical submanifold of $(N,c_N)$\,.\\
\indent
{\rm (ii)} If $\dim M\geq4$ then the following assertions are equivalent:\\
\indent
\quad{\rm (ii1)} $\Phi:(Q_M+Q,\J')\to(Q_N,\J'_N)$ is holomorphic.\\
\indent
\quad{\rm (ii2)} $\phi$ is geodesic.
\end{prop}
\begin{proof}
Let $x_0\in M$ and let $p_0\,, q_0$ be positive maximal isotropic spaces
which are tangent and normal, respectively, to $M$ at $x_0$. Let $S\subseteq M$ be a surface
such that $x_0\in S$ and one of the two isotropic directions tangent to $S$ at $x_0$
are contained in $p_0$\,; denote by $X_0$ a nonzero element of $T_{x_0}^{\C}\!S\cap p_0$
(obviously, $X$ is well-defined, up to some complex factor).\\
\indent
We may suppose that there exist two sections $p$ and $q$ of $Q_M$ and $Q$\,, respectively,
over $S$ which are horizontal at $x_0$ and such that $p_{x_0}=p_0$ and $q_{x_0}=q_0$\,.\\
\indent
{}From \cite[Proposition 2.6]{LouPan-II} it follows quickly that $\Phi:(Q_M+Q,\J)\to(Q_N,\J_N)$ is holomorphic
if and only if, for any $x_0\in M$ and any
such sections $p$ and $q$\,, we have $D^N_{\,X_0}Y\in p_0\oplus q_0$ and $D^N_{\,X_0}U\in p_0\oplus q_0$\,,
for any local sections $Y$ of $p$ and $U$ of $q$\,; equivalently, $c_N(D^N_{\,X_0}Y,U)=0$ for any local section
$Y$ of $p$ and $U$ of $q$\,. The proof of (i) follows quickly.\\
\indent
Similarly, $\Phi:(Q_M+Q,\J')\to(Q_N,\J'_N)$ is holomorphic if and only if, for any $x_0\in M$ and any
such sections $p$ and $q$\,, we have $D^N_{\,\overline{X_0}}Y\in p_0\oplus q_0$ and $D^N_{\,\overline{X_0}}U\in p_0\oplus q_0$\,,
for any local sections $Y$ of $p$ and $U$ of $q$\,. It follows that (ii1) is equivalent to
the fact that the Weyl connection induced by $D^N$ on $M$ is equal to $D^M$ and, for any $x_0\in M$
and any such sections $p$ and $q$\,, we have $c_N(D^N_{\,\overline{X_0}}Y,U)=0$ for any local sections
$Y$ of $p$ and $U$ of $q$\,. The proof of (ii) follows quickly.
\end{proof}

\indent
Similarly to the proof of Proposition \ref{prop:umbgeod}(ii)\,, we obtain the following:

\begin{rem}[cf.\ \cite{EelSal}\,]
If $\dim M=2$ then the equivalence (ii1)$\iff$(ii2)\,, of Proposition \ref{prop:umbgeod}\,, remains true if we replace (ii2)
with the following assertion:\\
\indent
(ii2) $M^2$ is a minimal surface in $(N,c_N,D^N)$\,.
\end{rem}

\section{On (1,1)-geodesic submanifolds} \label{section:(1,1)submfds}

\indent
Let $(N,c_N,D^N)$ be a Weyl space. For $1\leq r<\tfrac12\dim N$, let $\p_{N,r}:Q_{N,r}\to N$ be the bundle of
isotropic spaces on $(N,c_N)$ of complex dimension $r$\,. Denote by $\J_{N,r}$ and $\J'_{N,r}$ the almost CR-structures on
$Q_{N,r}$ whose eigenbundles corresponding to $-{\rm i}$ are constructed, similarly to $\tG$ and $\tG'$, respectively,
of Theorem \ref{thm:O'BrRaw}\,, by using the connection induced by $D^N$ on $Q_{N,r}$\,, and by taking $Q=N$ the trivial
bundle over $N$.\\
\indent
Note that, $\J_{N,r}$ does not depend of $D^N$ whilst $\J'_{N,r}$ determines $D^N$.
Furthermore, by Theorem \ref{thm:O'BrRaw}\,, the almost CR-structure $\J'_{N,r}$ is nonintegrable whilst,
if $r=1$ then $\J_{N,1}$ is integrable \cite{LeB-CR}\,. We shall prove the following result.

\begin{thm} \label{thm:LeBr-O'BrRaw}
The following assertions are equivalent, if $r\geq2$\,:\\
\indent
{\rm (i)} $\J_{N,r}$ is integrable.\\
\indent
{\rm (ii)} $(N,c_N)$ is flat.
\end{thm}
\begin{proof}
Assume $r\geq 2$ and let $R$ and $W$ be the curvature form of $D^N$ and the Weyl tensor of $(N,c_N)$\,, respectively.
We shall prove that the following assertions are equivalent:\\
\indent
\;\;(a) $R(\Lambda^2p)(p)\subseteq p$ for any $p\in Q_{N,r}$\,.\\
\indent
\;\;(b) $c_N(R(X,Y)X,Y)=0$ for any $X,Y\in T^{\C}\!N$ spanning an isotropic space.\\
\indent
\;\;(c) $W=0$\,.\\
Indeed, as any two-dimensional isotropic space on $(N,c_N)$ is contained in some $p\in Q_{N,r}$\,, we obviously have (a)$\Longrightarrow$(b)\,.
Also, (b)$\iff$(c) (see \cite{Pan-cf}\,) and, as $R(\Lambda^2p)(p)=W(\Lambda^2p)(p)$\,, for any isotropic space $p$ on $(N,c_N)$\,,
we have (c)$\Longrightarrow$(a)\,.\\
\indent
By Theorem \ref{thm:O'BrRaw}\,, we have (i)$\iff$(a)\,, and, by the Weyl theorem on flat conformal manifolds, (ii)$\iff$(c)\,.
The theorem is proved.
\end{proof}

\indent
Let $M\subseteq N$ be a submanifold, $\dim M=2r$\,. Let $c_M=c_N|_{_M}$ and let $D^M$ be the Weyl connection on $(M,c_M)$
induced by $D^N$. Also, let $\t'_M=(Q_M,M,\p_M,\J'_M)$ be the nonintegrable almost twistorial structure associated to $(M,c_M,D^M)$\,.\\
\indent
Suppose that there exists a section $p$ of $Q_{N,r}$ which is tangent to $M$.
As before, denote by $J^p$ the almost Hermitian structure on $(M,c_M)$ whose eigenbundle corresponding to $-{\rm i}$
is $p$\,.\\
\indent
Similarly to Proposition \ref{prop:2}\,, we obtain the following result (cf.\ \cite{Raw-ftwistor}\,).

\begin{prop} \label{prop:Weierstrass}
The following assertions are equivalent.\\
\indent
{\rm (i)} $p:(M,J^p)\to(Q_{N,r},\J'_{N,r})$ is holomorphic.\\
\indent
{\rm (ii)} $(M,J^p)$ is a $(1,1)$-geodesic submanifold of $(N,c_N,D^N)$ and the map $p:(M,J^p)\to(Q_M,\J'_M)$ is holomorphic.
\end{prop}

\begin{rem} \label{rem:Wass}
1) Proposition \ref{prop:Weierstrass} can be easily formulated in similar vein to Proposition \ref{prop:2}\,.\\
\indent
2) With the same notations as in Proposition \ref{prop:Weierstrass}\,, $p:(M,J^p)\to(Q_{N,r},\J_{N,r})$ is holomorphic
if and only if $J^p$ is integrable and $(M,J^p)$ is a $(2,0)$-geodesic submanifold of $(N,c_N,D^N)$\,.\\
\indent
3) In Proposition \ref{prop:Weierstrass}\,, assume that $(N,c_N,D^N)$ is the Euclidean space $\R^n$ with its canonical
conformal structure and flat connection. Then $Q_{N,r}=\R^n\times Q_{n,r}$ where $Q_{n,r}\subseteq{\rm Gr}_r(n,\C)$ is the manifold of
isotropic $r$-dimensional subspaces of $\C^{\!n}$.\\
\indent
Let $\widetilde{p}=\p_2\circ p:M\to F$ where $\p_2:\R^n\times Q_{n,r}\to Q_{n,r}$ is the projection. Then
$p:(M,J^p)\to(Q_{N,r},\J'_{N,r})$ is holomorphic
if and only if $\widetilde{p}:(M,J^p)\to Q_{n,r}$ is holomorphic.\\
\indent
Thus, by Proposition \ref{prop:Weierstrass}\,, \emph{$(M,J^p)$ is a $(1,1)$-geodesic submanifold of $(N,c_N,D^N)$ and
$p:(M,J^p)\to(Q_M,\J'_M)$ is holomorphic if and only if $\widetilde{p}:(M,J^p)\to Q_{n,r}$ is holomorphic}.
In the particular case $\dim M=2$\,, this gives $M^2$ minimal in $\R^n$ if and only if $\widetilde{p}$ holomorphic
which leads to the Weierstrass representation of minimal surfaces in Euclidean space.\\
\indent
4) A result similar to Proposition \ref{prop:umbgeod} can be easily written by working with the inclusion map
$Q_M\hookrightarrow Q_{N,r}$\,.
\end{rem}

\section{$f$-structures and pseudo horizontally conformal submersions} \label{section:f}

\indent
We start this section by recalling the following definition.

\begin{defn}[see \cite{BaiWoo2}\,,\,\cite{Bri}\,]
A map $\phi:(M,c)\to(N,J)$ from a conformal manifold to an almost complex manifold
is \emph{pseudo horizontally weakly conformal} if it pulls back $(1,0)$-forms on $N$
to isotropic $1$-forms on $(M,c)$\,. A map is \emph{pseudo horizontally conformal}
if it is submersive and pseudo horizontally weakly conformal.
\end{defn}

\begin{rem}
1) A submersion $\phi:(M,c)\to(N,J)$ from a conformal manifold to an almost complex manifold
is pseudo horizontally conformal if there exists an almost $f$-structure
$F$ on $(M,c)$ such that $T^0M={\rm ker}\dif\!\phi$ and $\phi:(M,F)\to(N,J)$ is holomorphic (cf.\ \cite{LouMo}\,).\\
\indent
2) Let $(M,c)$ be a conformal manifold and let $F$ be an almost $f$-structure on $M$. Then $F$ is an
$f$-structure on $(M,c)$ if and only if it is locally defined by pseudo horizontally conformal
submersions onto complex manifolds.
\end{rem}

\indent
Let $(M,c,D)$ be a Weyl space, $\dim M=m$\,. For $1\leq r<\tfrac12\,m$, let $\p_{M,r}:Q_{M,r}\to M$ be the bundle of
isotropic spaces on $(M,c)$ of complex dimension $r$\,. For $p\in Q_{M,r}$ let $F^p$ be the skew-adjoint
$f$-structure on $(T_{\p_{M,r}(p)}M,c_{\p_{M,r}(p)})$ whose eigenspace corresponding to $-{\rm i}$ is $p$\,.
Thus, $Q_{M,r}$ is also the bundle of skew-adjoint $f$-structures on $(M,c)$ with kernel of
dimension $m-2r$\,.\\
\indent
Let $\H$ be the connection induced by $D$ on $Q_{M,r}$ and let $T^0Q_{M,r}\subseteq\H$ be the subbundle characterised
by $\dif\!\p_{M,r}(T^0_pQ_{M,r})={\rm ker}F^p$, for all $p\in Q_{M,r}$\,. Also, let
$\tG_0\subseteq\H^{\C}$ be the subbundle such that $\dif\!\p_{M,r}\bigl((\tG_0)_p\bigr)$ is the eigenspace of $F^p$
corresponding to $-{\rm i}$\,, for all $p\in Q_{M,r}$\,. Denote $T^{0,1}Q_{M,r}=\tG_0\oplus({\rm ker}\dif\!\p_{M,r})^{0,1}$
and let $\Fa_{M,r}$ be the almost $f$-structure on $Q_{M,r}$ whose eigenbundles corresponding to $0$ and $-{\rm i}$ are
$T^0Q_{M,r}$ and $T^{0,1}Q_{M,r}$\,, respectively. Also, let $\Fa'_{M,r}$ be the almost $f$-structure on $Q_{M,r}$
whose eigenbundles corresponding to $0$ and $-{\rm i}$ are $T^0Q_{M,r}$ and $\tG_0\oplus({\rm ker}\dif\!\p_{M,r})^{1,0}$\,,
respectively.

\begin{rem}
1) Each of the almost $f$-structures $\Fa_{M,r}$ and $\Fa'_{M,r}$ determines $D$\,.\\
\indent
2) With the same notations as in Section \ref{section:(1,1)submfds}\,, the almost CR-structures induced
by $\Fa_{M,r}$ and $\Fa'_{M,r}$ are $\J_{M,r}$ and $\J'_{M,r}$\,, respectively.
\end{rem}

\indent
It is well-known (see \cite[Theorem 3.5]{Pan-Deva}\,) that if $m=3$ then $\Fa_{M,1}$ is integrable
if and only if $(M,c,D)$ is Einstein--Weyl. Also, from Theorem \ref{thm:O'BrRaw} it easily follows that $\Fa'_{M,r}$
is nonintegrable. We shall prove the following:

\begin{thm} \label{thm:tcc}
If $m\geq4$ then the following assertions are equivalent:\\
\indent
{\rm (i)} $\Fa_{M,r}$ is integrable.\\
\indent
{\rm (ii)} $D$ is, locally, the Levi-Civita connection of a constant curvature representative of $c$\,.
\end{thm}
\begin{proof}
Let $R$ be the curvature form of the connection induced by $D$ on $L^*\otimes TM$, where $L$ is the line bundle
of $M$. We claim that the following assertions are equivalent:\\
\indent
\;\;(a) $R\bigl(\Lambda^2(p^{\perp})\bigr)\bigl(p^{\perp}\bigr)\subseteq p^{\perp}$ for any $p\in Q_{M,r}$\,.\\
\indent
\;\;(b) $c(R(X,Y)X,Y)=0$ for any $X,Y\in T^{\C}\!M$ spanning a degenerate space.\\
\indent
\;\;(c) $(M,c,D)$ is flat and Einstein--Weyl.\\
Indeed, as any two-dimensional degenerate (and, if $r=1$\,, nonisotropic) space on $(M,c)$ is contained in $p^{\perp}$ for some $p\in Q_{M,r}$\,,
we obviously have (a)$\Longrightarrow$(b)\,. Also, by \cite{Pan-cf}\,, assertion (b) implies that $(M,c)$ is flat;
it follows quickly that (b)$\Longrightarrow$(c)\,.
By a result of M.~G.~Eastwood and K.~P.~Tod (\,\cite[Theorem 1]{EasTod-Crelle}\,; see \cite[Theorem 5.2]{Cal-F}\,), (c)$\iff$(ii)\,.
Clearly, (ii)$\Longrightarrow$(a) and the proof follows from Theorem \ref{thm:O'BrRaw}\,.
\end{proof}

\begin{rem}
By Theorems \ref{thm:LeBr-O'BrRaw} and \ref{thm:tcc}\,, if $\Fa_{M,r}$ is integrable then $\J_{M,r}$ is
integrable.
\end{rem}

\indent
Let $(M,g)$ be a Riemannian manifold of constant curvature such that $\Fa_{M,r}$ is simple. Then
there exists a holomorphic submersion from $(Q_{M,r},\Fa_{M,r})$ onto a complex manifold $Z_r(M)$ whose fibres
are the leaves of $T^0Q_{M,r}$\,. Then $Z_r(M)$ is the twistor space of $(Q_{M,r},M,\p_{M,r},\Fa_{M,r})$
(cf.\ \cite[\S6.8]{BaiWoo2}\,).

\begin{prop}[cf.\ \cite{LouPan-II}\,] \label{prop:tcc}
Let $p$ be a section of $Q_{M,r}$ and let $F^p$ be the corresponding almost $f$-structure on $(M,c)$\,.
The following assertions are equivalent:\\
\indent
{\rm (i)} $p:(M,F^p)\to(Q_{M,r},\Fa_{M,r})$ is holomorphic.\\
\indent
{\rm (ii)} $F^p$ is integrable and locally defined by pseudo horizontally conformal submersions with geodesic fibres
and for which the integrability tensor of the horizontal distribution is of degree $(1,1)$\,.
\end{prop}
\begin{proof}
Assertion (i) is equivalent to the fact that $D_XY\in\G(p^{\perp})$\,, for any $X,Y\in\G(p^{\perp})$\,; in particular,
if (i) holds then $F^p$ is integrable. Clearly, (i) is also equivalent to $D_XY\in\G(\overline{p}^{\,\perp})$\,, for any
$X,Y\in\G(\overline{p}^{\,\perp})$\,. Therefore, if (i) holds then
$p^{\perp}\cap\overline{p}^{\,\perp}\,(=(p\oplus\overline{p})^{\perp}={\rm ker}F^p)$ is geodesic.\\
\indent
Thus, if (i) holds then $F^p$ is integrable and locally defined by pseudo horizontally conformal submersions with geodesic fibres;
furthermore, if $X,Y\in\G(p)$ and $U\in\G\bigl((p\oplus\overline{p})^{\perp}\bigr)$ then, as $F^p$ is integrable, we have
$[U,X],[U,Y]\in\G(p^{\perp})$ and it follows that  $c(U,[X,Y])=-2c(D_UX,Y)=0$\,. This completes the proof of (i)$\Longrightarrow$(ii)\,.\\
\indent
By definition, $F^p$ integrable if and only if $p^{\perp}$ integrable. It follows that if $F^p$ is integrable then
$D_XY\in\G(p^{\perp})$\,, for any $X,Y\in\G(p)$\,. Also, if ${\rm ker}F^p\bigl(=(p\oplus\overline{p})^{\perp}\bigr)$ is geodesic then
$D_UV\in\G(p^{\perp})$\,, for any $U,V\in\G\bigl((p\oplus\overline{p})^{\perp}\bigr)$\,. Furthermore, an argument as above
shows that if $F^p$ is integrable then the integrability tensor of $(p\oplus\overline{p})^{\perp}$ is of degree $(1,1)$
if and only if $D_UX\in\G(p^{\perp})$\,, for any $X\in\G(p)$ and $U\in\G\bigl((p\oplus\overline{p})^{\perp}\bigr)$\,.
This completes the proof of (ii)$\Longrightarrow$(i)\,.
\end{proof}

\begin{rem}
Let $F$ be an $f$-structure on $M$. It is obvious that the almost CR-structure $T^{0,1}M$ is integrable if and only if
the integrability tensor of $T^{1,0}M\oplus T^{0,1}M$ is of degree $(1,1)$\,.
\end{rem}

\indent
{}From Proposition \ref{prop:tcc} we easily obtain the following result.

\begin{cor}[cf.\ \cite{LouPan-II}\,] \label{cor:hmtcc}
Let $p$ be a section of $Q_{M,1}$ and let $F^p$ be the corresponding almost $f$-structure on $(M,c)$\,.
The following assertions are equivalent:\\
\indent
{\rm (i)} $p:(M,F^p)\to(Q_{M,1},\Fa_{M,1})$ is holomorphic.\\
\indent
{\rm (ii)} $F^p$ is integrable and locally defined by submersive harmonic morphisms with geodesic fibres
(of codimension two).
\end{cor}

\indent
Let $(M,g)$ be a real analytic Riemannian manifold, $\dim M=m$\,. Then $(M,g)$ admits a (germ-unique) complexification $(M^{\C}\!,g^{\C})$\,.
Let $\p_{M^{\C}\!,r}:Q_{M^{\C}\!,r}\to M^{\C}$ be the bundle of $r$-dimensional isotropic spaces on $(M^{\C}\!,g^{\C})$\,.
If $1\leq r<\tfrac12\,m$\,, the complex version of Theorem \ref{thm:tcc} says that the following assertions are equivalent
(cf.\ \cite[\S2]{PanWoo-sd} and the references therein):\\
\indent
\;\;(i) For any $p\in Q_{M^{\C}\!,r}$ there exists a coisotropic and geodesic complex submanifold $S$ of $(M^{\C}\!,g^{\C})$ of (complex)
rank $m-2r$\,, with respect to $g^{\C}$\,, such that $T_{\p_{M^{\C}\!,r}(p)}S=p^{\perp}$.\\
\indent
\;\;(ii) $(M^{\C}\!,g^{\C})$ has constant (sectional) curvature.\\
\indent
Assume $(M^{\C}\!,g^{\C})$ (and, hence, also $(M,g)$\,) to be of constant curvature. Then, locally, the twistor space (in the sense
of \cite[Definition 2.1]{PanWoo-sd}\,) $Z_r(M^{\C})$ parametrizes the coisotropic geodesic (complex) submanifolds of $(M^{\C}\!,g^{\C})$ of rank $m-2r$\,.
It follows that, locally, we may assume $T^0Q_{M,r}$ simple and such that each of its leaves intersects the fibres of $\p_{M,r}$ at
most once (apply \cite[Remark 2.2(3)]{PanWoo-sd}\,). Then $Z_r(M)$ is an open submanifold of $Z_r(M^{\C})$\,; moreover,
$Z_r(M)$ is endowed with a holomorphic $m$-dimensional family of submanifolds each of which is holomorphically
diffeomorphic to the space of isotropic $r$-dimensional spaces on $\C^{\!m}$; the members of this family are called
the \emph{twistor submanifolds} of $Z_r(M)$ (see \cite[Remark 2.2(1)]{PanWoo-sd}\,).\\
\indent
We shall say that two submanifolds $S$ and $S'$ of a manifold $W$ are \emph{transversal} if $T_xS\cap T_xS'=\{0\}$\,, at each $x\in S\cap S'$.

\begin{cor} \label{cor:tcc}
Let $(M,g)$ be a Riemannian manifold of constant curvature and let $1\leq r<\tfrac12\,m$\,, where $m=\dim M$.\\
\indent
Then any pseudo horizontally conformal submersion, locally defined on $(M,g)$\,, with geodesic fibres of dimension $m-2r$
and for which the integrability tensor of the horizontal distribution is of degree $(1,1)$ corresponds, locally,
to a complex submanifold, of dimension $r$, of $Z_r(M)$ which is transversal to the twistor submanifolds.
\end{cor}
\begin{proof}
Any (local) pseudo horizontally conformal submersion $\phi$ on $(M,g)$ with connected geodesic fibres of dimension $m-2r$
and for which the integrability tensor of the horizontal distribution is of degree $(1,1)$ defines an $f$-structure
$F^{\phi}$ on $(M,g)$\,. Moreover, by Proposition \ref{prop:tcc}\,, $F^{\phi}$ corresponds to a holomorphic
section $p^{\phi}:(M,F^{\phi})\to(Q_{M,r},\Fa_{M,r})$\,. Hence, $T^0Q_{M,r}$ induces a foliation
on $p^{\phi}(M)$ whose leaves are mapped by $\p_{M,r}$ onto the fibres of $\phi$\,. Thus, locally, the projection
$Q_{M,r}\to Z_r(M)$ maps $p^{\phi}(M)$ onto a complex $r$-dimensional submanifold $N^{\phi}$ of $Z_r(M)$.
Then $\phi\longmapsto\!N^{\phi}$ gives the claimed correspondence.
\end{proof}

\begin{rem}
Let $(M,g)$ be a constant curvature Riemannian manifold and let $1\leq r<\tfrac12\,m$\,, where $m=\dim M$.\\
\indent
Then \emph{$Z_r(M)$ parametrizes naturally the set of pairs $(P,J)$ where $P$ is a totally geodesic submanifold of $(M,g)$\,, of codimension $2r$\,,
and $J$ is an orthogonal complex structure on the normal bundle of $P$ which is parallel with respect to the normal connection}.
(By \eqref{e:Ricci} and \eqref{e:Weyl}\,, the normal connection on the normal bundle of any totally umbilical submanifold of a
conformally-flat Riemannian manifold is flat.)\\
\indent
Let $\phi$ be a (local) pseudo horizontally conformal submersion on $(M,g)$ with connected geodesic fibres of dimension $m-2r$
and for which the integrability tensor of the horizontal distribution is of degree $(1,1)$\,. Let $N^{\phi}$ be the codomain of $\phi$
and let $J^{\phi}$ be the orthogonal complex structure
on $({\rm ker}\dif\!\phi)^{\perp}$ with respect to which $\dif\!\phi|_{({\rm ker}\dif\!\phi)^{\perp}}$ is holomorphic
at each point.\\
\indent
Then the correspondence of Corollary \ref{cor:tcc} is given by $\phi\longmapsto N^{\phi}$ where the inclusion map
$N^{\phi}\hookrightarrow Z_r(M)$ is defined by $y\mapsto\bigl(\phi^{-1}(y),J^{\phi}|_{\phi^{-1}(y)}\bigr)$\,, $(y\in N^{\phi})$.
\end{rem}

\indent
{}From Corollary \ref{cor:tcc} we obtain the following result of P.~Baird and J.~C.~Wood.

\begin{cor}[\,\cite{BaiWoo2}\,] \label{cor:BaiWoo}
Let $(M,g)$ be a Riemannian manifold of constant curvature. Then any submersive harmonic morphism, locally defined on $(M,g)$\,,
with geodesic fibres of codimension two corresponds, locally, to a complex one-dimensional submanifold of $Z_1(M)$
which is transversal to the twistor submanifolds.
\end{cor}

\indent
We end by describing the twistor spaces of the space forms $\R^m$, $S^m$ and $H^m$ (cf.\ \cite[\S6.8]{BaiWoo2}\,).
For this, we firstly describe the twistor spaces of the complex Euclidean space $\C^{\!m}$ and of the
complex unit hypersphere $S^{m}(\C)$\,.\\
\indent
Let $Q_{m,r}\subseteq{\rm Gr}_{m-r}(m,\C)$ be the space of coisotropic subspaces of $\C^{\!m}$ of rank $m-2r$\,.
We shall denote by the same symbol $Q_{m,r}$ its image through the complex analytic diffeomorphism
${\rm Gr}_{m-r}(m,\C)\to {\rm Gr}_r(m,\C)$ defined by $p\mapsto p^{\perp}$, for any  $p\in {\rm Gr}_{m-r}(m,\C)$\,.
Thus, $Q_{m,r}\subseteq{\rm Gr}_r(m,\C)$ is the space of isotropic subspaces of $\C^{\!m}$ of complex dimension $r$\,.
Let $E_{m,r}$ and $F_{m,r}$ be the restrictions to $Q_{m,r}$ of the tautological vector bundles on ${\rm Gr}_{m-r}(m,\C)$
and ${\rm Gr}_r(m,\C)$\,, respectively. As $Z_r(\C^{\!m})$ is the space of coisotropic planes in $\C^{\!m}$ of rank $m-2r$\,,
we have $Z_r(\C^{\!m})=\bigl(Q_{m,r}\times\C^{\!m}\bigr)/E_{m,r}=F_{m,r}^{\,*}$\,.\\
\indent
Similarly, $Z_r\bigl(S^m(\C)\bigr)$ is the space of (maximal) coisotropic geodesic submanifolds of $S^m(\C)$\,, of rank $m-2r$. As any such
submanifold is the intersection of $S^m(\C)$ with a coisotropic subspace, of rank $m-2r+1$\,, of $\C^{\!m+1}$,
we have $Z_r\bigl(S^m(\C)\bigr)=Q_{m+1,r}$\,.\\
\indent
It follows that $Z_r(\R^m)=F_{m,r}^{\,*}$, $Z_r(S^m)=Q_{m+1,r}$ and $Z_r(H^m)=Q_{m+1,r}\!\setminus\mathcal{C}_{m,r}$
for some closed set $\mathcal{C}_{m,r}\subseteq Q_{m+1,r}$\,.
To describe $\mathcal{C}_{m,r}$\,, consider the complex Euclidean space $\C^{\!m+1}$ as the complexification of the
Minkowski space $\R^{m+1}_1$ so that the complexification of $H^m\subseteq\R^{m+1}_1$
to be the complex hypersphere, of radius the imaginary unit. Then $\mathcal{C}_{m,r}$ is the set of coisotropic subspaces
$p\subseteq\C^{\!m+1}$ of rank $m-2r+1$ such that $p^{\perp}\cap\R^{m+1}_1\neq\{0\}$\,.

\end{document}